\documentstyle[12pt,amscd]{amsart}
\newfont{\ams}{msbm10 at 12pt}
\newfont{\amsi}{msbm8}
\newcommand{\fu}{{\frak u} }

\newcommand{\cC}{{\cal C} }

\newcommand{\BZ}{{\Bbb Z} }
\newcommand{\BC}{{\Bbb C} }
\newcommand{\hL}{\hat{L}}
\newcommand{\hP}{\hat{P}}
\newcommand{\ad}{\hbox{\bf{ad}}}

\newcommand{\Hom}{\operatorname{Hom}}
\newcommand{\Ext}{\operatorname{Ext}}
\newcommand{\vareps}{\varepsilon}
\newcommand{\Image}{\operatorname{Im}}
\newcommand{\Ker}{\operatorname{Ker}}
\renewcommand{\mod}{\operatorname{mod}\ }
\begin{document}
\title{Decomposition of the adjoint representation
of the small quantum $sl_2.$}
\author{Victor Ostrik}
\address{Independent University of Moscow, PO Box 230,
Moscow, 117463, Russia}
\date{December 1995}
\maketitle
\section{Introduction.}

     \subsection{}  Given  a  finite  type  root  datum  and  a
primitive root of unity $q=\sqrt[l]{1}$, G.~Lusztig has defined
in [Lu] a remarkable finite dimensional Hopf algebra $\fu$ over
the  cyclotomic  field ${\Bbb Q}(\sqrt[l]{1})$. It is called  a
{\em restricted  quantum  universal  enveloping  algebra}, or a
{\em small quantum group}.

     Recall that for a Hopf algebra $A$ with a coproduct
$\Delta,$ antipode $S,$ and counit $\vareps$ the adjoint
representation is defined in the following way. $A$ is an
$A$-bimodule with respect to left and right multiplication.
Using the antipode we can consider $A$ as $A\otimes A$-module.
Combining this with the coproduct we get a new structure of
$A$-module on $A.$ It is called the {\em adjoint
representation} and denoted by $\ad$.

     \subsection{} In this note we study the adjoint
representation of $\fu$ in the simplest case of the root datum
$sl_2$.

     The semisimple part of this representation is of big
importance in the study of local systems of conformal blocks in
WZW model for $\hat{sl}_2$ at level $l-2$ in arbitrary genus.
The problem of distinguishing the semisimple part is closely
related to the problem of integral representation of conformal
blocks (see [BFS]).

     We find all the indecomposable direct summands of $\ad$
with multiplicities. To formulate the answer let us recall a
few notations from the representation theory of $\fu$. The
representation $\ad$ is naturally ${\Bbb Z}$-graded, so we
consider the category $\cC$ of ${\Bbb Z}$-graded $\fu$-modules.
The simple modules in this category are parametrized by their
highest weights which can assume arbitrary integer values. The
simple module with a highest weight $\lambda\in{\Bbb Z}$ is
denoted by $L(\lambda)$, and its indecomposable projective
cover is denoted by $P(\lambda)$.

     \subsection{} {\bf Main theorem.} {\em Let $l\geq 3$ be an
odd integer.
     The adjoint representation $\ad$ is isomorphic to the
direct sum of the
     modules $P(0),P(2),\ldots,P(l-3);
L(-l+1),L(-l+3),\ldots,L(2l-4), L(2l-2)$ with the following
multiplicities:

     (i) the multiplicity of $P(l-1)\simeq L(l-1)$ is $l;$

     Let $i\in [0,\ldots,\frac{l-3}{2}].$

     (ii) the multiplicity of $P(2i)$ is $\frac{l+1}{2}+i;$

     (iii) the multiplicity of $L(2i)$ is $l-1-2i;$

     (iiii) the multiplicities of $L(2l-2-2i)$ and $L(-2-2i)$
are $\frac{l-1}{2}-i.$}

     \subsection{} G.~Lusztig has defined (see e.g. [Lu]) a
{\em quantum enveloping algebra with divided powers} $U$
containing $\fu$ as a Hopf subalgebra. Its finite dimensional
irreducible representations are parametrized by their highest
weights which can assume arbitrary nonnegative integer values.
The simple module with a highest weight $\lambda\in{\Bbb N}$ is
denoted by $\hL(\lambda)$, and its indecomposable projective
cover is denoted by $\hP(\lambda)$. There is a natural
restriction functor from the category of finite dimensional
$U$-modules to $\cC$. It appears that $\ad$ lies in its
essential image.

     {\bf Corollary.} {\em The structure of $\BZ$-graded
$\fu$-module on $\ad$ can be lifted to the structure of
$U$-module isomorphic to the direct sum of the modules
$\hP(0),\hP(2),\ldots,\hP(l-3); \hL(0),\hL(2),\ldots,\hL(2l-4),
\hL(2l-2)$ with the following multiplicities:

     (i) the multiplicity of $\hP(l-1)\simeq \hL(l-1)$ is $l;$

     Let $i\in [0,\ldots,\frac{l-3}{2}].$

     (ii) the multiplicity of $\hP(2i)$ is $\frac{l+1}{2}+i;$

     (iii) the multiplicity of $\hL(2i)$ is $l-1-2i;$

     (iiii) the multiplicity of $\hL(2l-2-2i)$ is
$\frac{l-1}{2}-i.$}

     \subsection{} The author is very grateful to Michael
Finkelberg; it is hard to overestimate
     his help with writing down the present work. Also I am
happy to express my sincere gratitude to Thomas Kerler who
taught me everything I know about $\ad$.

     \section{Notations}

     \subsection{} In this section we recall the necessary
facts about $\fu$ and its representations, following mainly
[Lu].

     Let $l$ be an odd integer, $l>1.$ Let $q\in {\Bbb C}$ be a
primitive $l-$th root of unity. Let
$(i)_q=\frac{q^i-q^{-i}}{q-q^{-1}}.$

     Let $\fu$ be an associative algebra over ${\Bbb C}$ with
generators $E, F, K, K^{-1}$ and relations:

     $$ KK^{-1}=K^{-1}K=1; $$

     $$ KEK^{-1}=q^2E, KFK^{-1}=q^{-2}F; $$

     $$ EF-FE=\frac{K-K^{-1}}{q-q^{-1}};$$ $$ E^{l}=F^{l}=0,
K^{l}=1.$$

     The algebra $\fu$ is finite-dimensional and $\dim\fu=l^3.$

     Let $\omega$ be an automorphism of the algebra $\fu$ given
on generators by the formulas: $$
     \omega (E)=F,\ \omega (F)=E,\ \omega (K)=K^{-1}. $$

     The algebra $\fu$ is Hopf algebra with respect to
coproduct $\Delta,$ antipode $S$ and counit $\varepsilon$ given
by the formulas: $$\Delta (E)=E\otimes 1+K\otimes E, \Delta
(F)=F\otimes K^{-1}+1\otimes F, \Delta (K)=K\otimes K;$$ $$
S(E)=-K^{-1}E, S(F)=-FK, S(K)=K^{-1};$$ $$ \varepsilon
(E)=\varepsilon (F)=0, \varepsilon (K)=1.$$

     \subsection{} \label{category} Let $\cC$ be a category of
finite-dimensional $\BZ$-graded $\fu$-modules
$V=\bigoplus\limits_{i\in \BZ} V^i$ such that the following
conditions hold:

     (a) $E$ is operator of degree 2, i.e. $E$ acts from $V^i$
to $V^{i+2};$

     (b) $F$ is operator of degree $-2$;

     (c) $K$ acts on $V^i$ by multiplication by $q^i.$

     The morphisms in category $\cC$ are morphisms of
$\fu$-modules compatible with $\BZ$-grading.

     \subsection{} We introduce the duality $D$ on category
$\cC.$ If
     $V\in \cC ,$ then $D(V)$ is $V^*$ as a vector space. The
action of $x\in \fu$ on $D(V)$ is given by the formula
$(xf)(v)=f(\omega S(x)v),$ where $f\in D(V), v\in V, S$ is the
antipode.

     \subsection{} \label{adjoint} Let us define the adjoint
representation $\ad \in \cC$ (see e.g. [LM]).
     Let $x$ be an element of $\fu.$ The adjoint action of
generators is given by the folloing formulas: $$
ad(E)x=Ex-KxK^{-1}E=K[K^{-1}E,x], $$ $$ ad(F)x=FxK-xFK=[F,x]K,
$$ $$ ad(K)x=KxK^{-1}. $$

     \subsection{} \label{grading}
     Now we introduce $\BZ$-grading on adjoint representation.
We put $\deg(E)=2, \deg(F)=-2, \deg(K)=0$ and
$\deg(ab)=\deg(a)+\deg(b)$ for any $a,b\in \fu $ such that
$\deg(a), \deg(b)$ are defined. Note that all the weights of
$\ad$ are even integers in the interval $[2-2l,\ldots,2l-2].$

     \subsection{} \label{autoduality}
     It is known (see e.g. [LM] or [BFS]) that $D(\ad)\simeq
\ad .$

     \section{$\fu$-modules.} \subsection{}
     It is easy to check that an element \begin{equation}
\label{Casimir}
X=EF+\frac{q^{-1}K+qK^{-1}}{(q-q^{-1})^2}=FE+\frac{qK+q^{-1}K^{
-1}}{(q-q^{-1})^2}
     \end{equation} lies in the center of algebra $\fu$ (see
e.g. [Ke]). The element $X$ is called Casimir element. It
satisfies the following equation of degree $l$ (see loc.cit.):
\begin{equation} \label{Cequat}
     P(X):=\prod_{j\in \BZ/l\BZ} (X-b_j)=0, \end{equation}
where $b_j=\frac{q^{j+1}+q^{-j-1}}{(q-q^{-1})^2}$ ($q^l=1,$ so
$q^j$ is well defined). In particular $b_j=b_{j'}$ if
$j+j'=l-2.$ The root $b_{-1}=\frac{2}{(q-q^{-1})^2}$ of $P$ has
multiplicity 1, and the rest roots $b_j$ have multiplicity 2.

     \subsection{} Let $j\in \BZ/l\BZ.$ Let $\cC_j$ be a full
subcategory of $\cC$ such that $(X-b_j)$ acts nilpotently on
objects of $\cC_j.$ In what follows we will identify $\cC_j$
and $\cC_{j'}$ if $j+j'=l-2.$

     \subsection{} Let us fix $H'$--- a maximal subset of
$\BZ/l\BZ$ with the following property: if $\{j,j'\}$ is
two-element subset of $H'$ then $j+j'\ne l-2.$ We have
$H'=\{-1\}\cup H,$ where $H=H'-\{-1\}.$

     \subsection{} For any $j\in H$ we define the integers
$J,J'$ by the following
     properties:

     (1) $0\le J<J'<l;$

     (2) $J+J'=l-2;$

     (3) $(J-j)(J'-j)\equiv 0~(\mod l).$

     \subsection{} The category $\cC$ is a direct sum of
subcategories $\cC_j$ where $j$ runs through $H'.$

     \subsection{} \label{Verma}
     We denote by $\fu^{\pm}\subset \fu$ the subalgebra
generated by $K,E$ (resp. $K,F$). For $\lambda \in \BZ$ denote
by $\BC_{\lambda}^{\pm}$ the one-dimensional $\BZ-$graded
$\fu^{\pm}-$module of weight $\lambda$ such that $K$ acts as
$q^{\lambda}$ and $E$ (resp. $F$) acts as zero on it. We denote
by $M^{\mp}(\lambda)$ the $\BZ-$graded $\fu-$module $\fu
\otimes_{\fu^{\pm}}\BC_{\lambda}^{\pm}.$ The modules
$M^{\pm}(\lambda)$ are called Verma modules. Let
$M^{\pm}(\lambda)\ni v^{\pm}(\lambda):=1\otimes 1.$ Let $V\in
\cC, v\in V^{\lambda}$ and $E\cdot v=0$ (resp. $Fv=0$). Then
$v$ is called an upper singular (resp. a lower singular) vector
in $V,$ and there exists a unique morphism $\phi :
M^{\pm}(\lambda)\to V$ such that $\phi (v^{\pm}(\lambda))=v.$

     \subsection{} \label{modules}
     For each $\lambda \in \BZ$ there is a unique up to
isomorphism simple module $L(\lambda )\in \cC $ with highest
weight $\lambda .$ The modules $L(\lambda_1)$ and
$L(\lambda_2)$ are isomorphic iff $\lambda_1\equiv \lambda_2
(\mod l).$ Indecomposable projective cover of $L(\lambda)$ will
be denoted by $P(\lambda).$
     We have $D(L(\lambda))\simeq L(\lambda)$ and
$D(P(\lambda))\simeq P(\lambda).$ In particular $P(\lambda)$ is
injective; $\dim \Hom(L(\lambda),P(\lambda))=1.$

     \subsection{} \label{Steinberg}
     The set of isomorphism classes of simple objects in
category
     $\cC_{-1}$ is $\{L(\lambda), \lambda \equiv -1(\mod l)\}.$
As $\fu$-modules without grading all the $L(\lambda)$ are
isomorphic to one and the same $\fu$-module St (Steinberg
module). It has dimension $l.$ The category $\cC_{-1}$ is
semisimple. In particular $P(\lambda)=L(\lambda).$

     \subsection{} \label{nonsem} Let $j\in H.$
     The set of isomorphism classes of simple modules in
$\cC_j$ is $\{L(\lambda), \lambda \equiv J(\mod l)$ or $\lambda
\equiv J'(\mod l)\}.$ The modules $L(\lambda), \lambda \equiv
J(\mod l)$
     (resp. $\lambda \equiv J'(\mod l)$) are isomorphic as
$\fu$-modules. Their dimension is $J+1$ (resp. $J'+1$).

     The projective module $P(\lambda)$ admits a filtration
$P(\lambda)\supset W(\lambda)\supset L(\lambda)\supset 0$ such
that $P(\lambda)/W(\lambda)\simeq L(\lambda),
W(\lambda)/L(\lambda)\simeq L(\lambda')\oplus L(\lambda'')$
where $\lambda' \ne \lambda'',$ $\lambda'\equiv \lambda''\equiv
-2-\lambda ~(\mod l), |\lambda-\lambda'|<2l>
|\lambda-\lambda''|.$
     In particular $\dim P(\lambda)=2l.$

     \subsection{} It is easy to see from
{}~\ref{grading},~\ref{Steinberg} and ~\ref{nonsem}
     that all the simple subquotients of adjoint representation
have the type $L(\lambda)$ where $\lambda$ is an even integer
from the interval $[1-l,\ldots,2l-2].$ Hence a projective
module $P(\lambda)$ can be a subquotient of $\ad$ only if
$\lambda \in 2\BZ \cap [0,\ldots,l-1].$ In particular each
subcategory $\cC_j$ contains only one isomorphism class of such
projectives.

     \subsection{} \label{lemma} {\bf Lemma.} {\em Suppose
$V\in \cC$ is indecomposable and the action of Casimir element
$X$ on $V$ is not semisimple. Then there exists $\lambda
\not\equiv -1~(\mod l)$ such that $V\simeq P(\lambda).$}

     {\bf Proof.} Casimir acts nonsemisimply on regular
representation (see ~(\ref{Cequat})). It follows that action on
projective modules $P(\lambda), \lambda \not\equiv -1~(\mod l)$
is not semisimple. It is easy to see that the space of
eigenvectors of Casimir in $P(\lambda)$ is $W(\lambda).$

     Let $W\subset V$ be a maximal submodule of $V$ such that
$X$ acts on $W$ semisimply. Choose $0\ne \varphi \in \Hom
(V,L)$ where $L=L(\lambda)$ for some $\lambda \in \BZ$ such
that Ker $\varphi$ contains $W.$ We have a morphism $\psi \in
\Hom (P,V)$ where $P=P(\lambda)$
     such that the diagram is commutative:

     \begin{equation*} \begin{CD}
     P \\ @V\operatorname{\psi}VV \\
     V @>\varphi>> L @>>> 0 \\ \end{CD} \begin{picture}(0,0)
     \put(-93,14){\vector(4,-3){32}} \end{picture}
     \end{equation*}

     If Ker $\psi \ne 0$ then X acts on Im $\psi$ semisimply.
Therefore $W$ is not maximal. We have a contradiction. If Ker
$\psi =0$ then we have injection $P\hookrightarrow V.$ From
{}~\ref{modules} follows that $P$ is direct summand of $V.$ The
proof is complete. $\Box$

     \section{The blocks of adjoint representation}
\label{grad}
     In what follows we always will identify $S_j$ and $S_{j'}$
where $j+j'=l-2$ and $S$ is an object, map, etc. Also we will
identify $S_J$ and $S_j$ if $J\in \BZ, J\equiv j~(\mod l).$

     \subsection{} The regular action of Casimir X (by
multiplication) is an endomorphism of adjoint
     representation. This gives a decomposition of adjoint
representation into {\em blocks} $\ad =\bigoplus\limits_{j\in
H'}\ad _j$ where $(X-b_j)$ acts nilpotently on $\ad _j.$ Let
$pr_j$ denote a projection onto $\ad_j.$

     \subsection{} Let $j\in H.$ Let $M_j=\Ker (X-b_j),$
$N_j=\ad_j\cap \Image (X-b_j).$
     Each $\ad_j$ admits a filtration $\ad_j\supset M_j\supset
N_j\supset 0.$ The rest of this section is a computation of
assotiated graded of this filtration. Evidently
$\ad_j/M_j\simeq N_j.$ It remains to compute $N_j, M_j/N_j.$
     It is convenient to put $N_{-1}=\ad_{-1}.$

     \subsection{} Recall (see ~\ref{category}) that
$\ad^0\subset \ad$ denotes the zero weight
     space. Let $\ad_j^0=\ad_j\cap \ad^0$ (for all $j\in H'$),
$N_j^0=N_j\cap \ad^0, M_j^0=M_j\cap \ad^0$
     (for $j\in H$) and $N_{-1}^0= \ad_{-1}^0.$ We will compute
the action of ad$(X)$ on $\ad_j^0.$ We start with a computation
of action of ad$(X)$ on $N_j^0.$

     \subsection{} \label{dimensions} {\bf Lemma.} {\em We have

     (a) $\dim \ad_j=2l^2$ if $j\in H$ and $\dim \ad_{-1}=l^2;$

     (b) if $j\in H$ then $\dim N_j=\dim
\ad_j/M_j=(J+1)^2+(J'+1)^2$ and $\dim M_j/N_j=4(J+1)(J'+1).$

     (c) $\dim \ad^{2m}=l(l-|m|)$ for all $m\in \BZ$ such that
$|m|<l;$

     (d) $\dim \ad_j^0=2l$ for $j\in H$ and $\dim N_j^0=l$ for
all $j\in H'.$

     (e) $\dim \ad_j^{2m}\ge 2(l-|m|)$ if $j\in H$ and $|m|\ge
l-1-J.$ }

     {\bf Proof.} (a), (b), (c), (d) are trivial. Let us prove
(e). Suppose $m>0.$ It is easy to see from consideration of
$\fu$-action on Verma modules $M^+(J)$ and $M^+(J')$ that $E^m$
acts nontrivially
     at least on $2(l-m)$ weights. By standard arguments with
Vandermond determinant we obtain that the desired dimension is
at least $2(l-m).$
     The proof for $m<0$ is similar. $\Box$

     \subsection{} \label{freedom}
     The subspace $\ad^0\subset \fu$ is a subalgebra of $\fu.$
It is generated as
     algebra by $K$ and $X$ (see [Ke]). Moreover $\ad^0$ is a
free module over a subalgebra generated by $K$ (see loc. cit.).
In particular we have $M_j^0=N_j^0.$

     We have \begin{equation} \label{Casact}
\mbox{ad}(X)K^i=\frac{q^{2i-1}+q^{1-2i}}{(q-q^{-1})^2}K^i-(q^i-
q^{-i})^2XK^{i+1}
     +(i)_q(i+1)_qK^{i+2} \end{equation}

     \subsection{} \label{mat1}
     For $j\in H$ the elements $ (X-b_j)pr_j K^i, i=1,\ldots,l
$ (resp. $ pr_{-1} K^i, i=1,\ldots,l $ for $j=-1$)
     form a basis of $N_j^0.$ In this basis ad$(X)$ acts as a
lower-triangular matrix:
     \begin{equation} \label{matrix} A(j)=\left(
\begin{array}{ccccc} b_0 & 0 & 0 &\ldots & 0 \\
     (q-q^{-1})^2b_j & b_2& 0 &\ldots &0 \\ (1)_q(2)_q &
(q^2-q^{-2})^2b_j& b_4 &\ldots & 0 \\ 0 & (2)_q(3)_q &
(q^3-q^{-3})^2b_j &\ldots & 0 \\ \vdots &\vdots & \vdots
&\ddots &\vdots \\ 0 & 0 & 0 &\ldots & b_0 \\ \end{array}
\right) \end{equation}

     \subsubsection{} \label{basis} {\em Remark.} The vectors
$pr_jK^i,i=1,\ldots,l; (X-b_j)pr_jK^i,
     i=l+1,\ldots,2l$ form a basis of $\ad_j^0$ ($j\in H$).

     \subsection{} Let $k\in 2\BZ \cap [0,\ldots,l-1].$ The
eigenvalues of this matrix are $b_k$ (with multiplicity 2 if
$k\ne l-1$
     and multiplicity 1 if $k=l-1$). Let $k\ne l-1.$ It is
obvious that there exists 1 or 2 eigenvectors of $A(j)$
corresponding to eigenvalue $b_k.$ We have 2 eigenvectors iff
the determinant $d(j,k)$ of matrix (see Lemma ~\ref{lin1})
\begin{equation} D(j,k)=\left( \begin{array}{cccc}
(q^{k/2+1}-q^{-k/2-1})^2b_j & b_{k+2}-b_k & 0 & \ldots \\
(k/2+1)_q(k/2+2)_q &(q^{k/2+2}-q^{-k/2-2})^2b_j & b_{k+4}-b_k &
\ldots \\ \vdots &\vdots &\ddots &\vdots \\ 0 &\ldots &\ldots
&\ldots \\ \end{array} \right) \end{equation}
     is equal to zero. It is easy to see that this determinant
is a polynomial in $b_j^2$ of degree $\frac{l-1-k}{2}.$
     Since $b_j^2=b_{j'}^2 \Rightarrow b_j=b_{j'},$ the
polynomial $d(j,k)$ vanishes for at most $\frac{l-1-k}{2}$
values of $j\in H'.$

     \subsection{} \label{halyava}
     In this section we prove the following proposition:

     {\bf Proposition.} {\em For any $j\in H$ we have the
following decomposition: $$ N_j\simeq \bigoplus_{i=0}^{J}
(L(2i)\oplus L(2i))
     \oplus \bigoplus_{i=J+1}^{\frac{l-3}{2}}P(2i) \oplus
L(l-1). $$ }

     {\bf Proof.} The proof proceeds by induction: we start
from $j=\frac{l-3}{2},$ then proceed to $j=\frac{l-5}{2}$ etc.

     It follows from ~\ref{mat1} that for any $j\in H'$ the
module $N_j$ contains as subquotients $L(0),L(2),\ldots,
     L(l-3)$ with multiplicities 2 and $L(l-1)$ with
multiplicity 1 (since only these modules have nontrivial zero
weight space).

     {\bf Lemma.} {\em Let $j+1=\frac{l-1}{2}.$ Then $N_j\simeq
L(0)\oplus L(0) \oplus \ldots \oplus L(l-3)\oplus L(l-1).$}

     {\bf Proof.} Let us compute the dimensions. We have
     $\dim
N_j=(\frac{l-1}{2})^2+(\frac{l+1}{2})^2=\frac{l^2+1}{2}$ (see
Lemma ~\ref{dimensions}(b)). On the other hand by the above
     we have $\dim N_j\ge 2\dim L(0)+2\dim L(2)+\ldots +2\dim
L(l-3)+ \dim L(l-1)=2\cdot 1+2\cdot 3+\ldots +2\cdot
(l-2)+l=\frac{l^2+1}{2}.$
     It follows that in this case $N_j$ is a direct sum of
$L(0),L(2),\ldots, L(l-3)$ with multiplicities 2 and $L(l-1)$
with multiplicity 1 (since $\Ext^1(L(\lambda),L(\mu))=0 \quad
\forall \lambda,\mu \in \{0,2,\ldots,l-1\}$). $\Box$

     The Lemma implies that all eigenvalues of
$A(\frac{l-3}{2})$ are semisimple. It follows from
     ~\ref{mat1} that eigenvalue $b_{l-3}$ is not semisimple
for all the rest $j.$ Therefore for $j\ne \frac{l-3}{2}$ the
corresponding $N_j$ contains projective
     submodule $P(l-3)$ (see ~\ref{lemma}).

     Now let $j+1=\frac{l-3}{2}.$ Then $\dim
N_j=\frac{l^2+9}{2}.$ On the other hand $\dim N_j\ge 2\cdot
1+\ldots +2\cdot (l-4)+2l+l=\frac{l^2+9}{2}.$ So in this case
$N_j$ is a direct sum of $L(0),L(2),\ldots, L(l-5)$ with
multiplicities 2, $L(l-1)$ with multiplicity 1, and $P(l-3).$
     As above it follows that all the rest $N_j$ contains
projective submodules $P(l-5)$ and $P(l-3)$ etc. The
Proposition is proved. $\Box$

     \subsubsection{} \label{block} {\bf Corollary.} {\em We
have: $$ \ad_{-1}=N_{-1}=
\bigoplus\limits_{i=0}^{\frac{l-1}{2}} P(2i). $$}

     {\bf Proof.} It follows from the proof of the Proposition
{}~\ref{halyava} that all the eigenvalues of the matrix $A(-1),$
except for $b_{l-1},$ are not semisimple. Hence the result
follows from the Lemma ~\ref{lemma} and computation of
dimensions. $\Box$

     \subsection{} \label{forrem}
     The Corollary ~\ref{block} gives a decomposition of
$\ad_{-1}.$ So in what follows we will assume that $j\ne -1$
i.e. $j\in H.$

     {\bf Lemma.} {\em The module $M_j/N_j$ is a direct sum of
modules $L(\lambda)$ with multiplicity 2 where $\lambda$ is
even and satisfies one of the following conditions:

     (i) $\lambda$ lies in the interval
$[2(l-J+1),\ldots,2l-2];$

     (ii) $\lambda$ lies in the interval $[-2J-2,\ldots,-2].$}

     {\bf Proof.} Recall that $M_j^0=N_j^0.$ Hence $M_j/N_j$
     contains only subquotients $L(\lambda)$ where either
$\lambda >l$ or $\lambda <0.$ By the Lemma
{}~\ref{dimensions}(e), Proposition ~\ref{halyava},
     and Corollary ~\ref{block}, we have $\dim \ad_j^{2m}\ge
2(l-|m|)$ if $j\in H$ and $\dim \ad_{-1}^{2m}\ge l-|m|$
     for any $m\ge \frac{l+1}{2}.$ It follows from the Lemma
{}~\ref{dimensions}(c) that $\dim \ad_j^{2m}=2(l-|m|)$ for any
$j\in H$ and $m\ge \frac{l+1}{2}.$

     Hence for any $m\ge \frac{l+1}{2}$ we have exactly two
upper (resp. lower) singular vectors of weight $2m$ (resp.
$-2m$). It follows that $\ad_j$ has two simple subquotients
with highest weight $2m$ and two simple subquotients with
lowest weight $-2m$ for any $m\ge \frac{l+1}{2}.$
     Thus $\ad_j$ has the following subquotients: $L(2i)$ where
$i\in [0,\ldots,\frac{l-3}{2}]$ with multiplicities 4; $L(l-1)$
with multiplicity 2; $L(-2-2i)$ and $L(2l-2-2i)$, where $i\in
[0,\ldots,\frac{l-3}{2}]$, with multiplicities 2. The
computation of dimensions shows that these modules are all the
subquotients of $\ad_j.$
     It follows from Proposition ~\ref{halyava} that
$[N_j:L(\lambda)]=1$ if $\lambda =-2-2i$ and $\lambda
=2l-2-2i$, where $i \in [J+1,\ldots,\frac{l-3}{2}].$ Since
$\ad_j/M_j\simeq N_j,$ any simple subquotient of $M_j/N_j$ is
of the type $L(\lambda)$ where $\lambda \in
[2(l-J+1),\ldots,2l-2] \cup [-2J-2,\ldots,-2]$ is even. But for
any $\lambda ,\mu $ satisfying such conditions we have
$\Ext^1(L(\lambda), L(\mu))=0.$ The result follows. $\Box$

     \subsubsection{} {\em Remark.} It follows from the proof
of the Lemma ~\ref{forrem}
     that $\dim \ad_j^{2m}=2(l-|m|)$ for any $j\in H$ and $m\in
\BZ, |m|<l.$

     \subsection{} \label{resume}
     Let $k$ be an even integer and $k\in [0,l-1].$ Let
$\ad_j(k)$ be a summand corresponding to the subcategory
$\cC_k$ in $\ad_j.$ Let $M_j(k)=M_j\cap \ad_j(k)$ and
$N_j(k)=N_j\cap \ad_j(k).$ Let us summarize the results of the
present section.

     \subsubsection{} $\ad_{-1}(k)=P(k).$

     \subsubsection{} \label{resume1}
     If $j\in H$ then

     (a) $\ad_j(l-1)=L(l-1)\oplus L(l-1);$

     (b) if $k\ge 2J+2$ then $\ad_j(k)=P(k)\oplus P(k);$

     (c) if $k\le 2J$ then $\ad_j(k)$ admits a fitration
$\ad_j(k)\supset M_j(k) \supset N_j(k)\supset 0$ with the
following associated graded factors: $$
     N_j(k)\simeq L(k)\oplus L(k); $$ $$
     M_j(k)/N_j(k)\simeq L(2l-2-k)\oplus L(2l-2-k)\oplus
L(-2-k)\oplus L(-2-k); $$ $$
     \ad_j(k)/M_j(k)\simeq L(k)\oplus L(k). $$

     \section{The proof of the main theorem} \subsection{}
\label{mat2}
     Let us find the multiplicities of projective submodules in
$\ad_j.$ Recall (see Remark ~\ref{basis}) that the vectors
     $pr_jK^i,i=1,\ldots,l; (X-b_j)pr_jK^i, i=l+1,\ldots,2l$
form a
     basis of $\ad_j^0.$ In this basis ad$(X)$ acts as a block
matrix $$ A'(j)=\left( \begin{array}{cc}
     A(j) & 0 \\ B & A(j) \\ \end{array} \right) $$ where
$A(j)$ is a matrix ~(\ref{matrix}) and $B$ is a matrix (see
{}~(\ref{Casact})) $$ \left( \begin{array}{cccc} 0 & 0 &0 &
\ldots \\ -(q-q^{-1})^2 & 0 &0 & \ldots \\ 0 & -(q^2-q^{-2})^2
& 0 & \ldots \\ \ldots & \ldots & \ldots & \ldots \\
\end{array} \right) $$

     \subsection{} Let $b_k$ be a nonsemisimple eigenvalue of
$A(j).$ Then a summand corresponding to the subcategory $\cC_k$
in $\ad_j$ is a sum of 2 copies of projective $P(k).$

     \subsection{} Let $b_k$ ($k\ne -1$) be a semisimple
eigenvalue of $A(j),$ i.e. $k\le 2J.$

     \subsubsection{} {\bf Lemma.} {\em In this case $\ad_j$
contains projective module from
     category $\cC_k.$}

     {\bf Proof.} It is enough to prove that the matrix $A'(j)$
has exactly 3 eigenvectors with eigenvalue $b_k$ or,
equivalently, that the matrix $A'(j)-b_k$ has corank 3. Let us
denote by $\tilde A'(j,k)$ the matrix $A'(j,k)-b_k$ with the
$i-$th and $(i+l)-$th columns divided by positive numbers
     $-(q^i-q^{-i})^2$ for any $i \in [1,\ldots,l-1].$ In order
to apply the Lemma ~\ref{lin2} let us put $A'=\tilde A'(j,k).$
Then corresponding matrix $D$ in notations of the Lemma
{}~\ref{lin2} is the matrix $D(j,k)$ with columns divided by some
positive numbers. In order to check the semisimplicity of the
matrix $D$ put $q=\exp (\pi i\frac{l+1}{l}).$ Then conditions
of Lemma ~\ref{lin3} hold. Indeed, off-diagonal entries of
$D(j,k)$ are either $(t)_q(t+1)_q$ or $b_{k+2t}-b_k=
(t)_q(t+k+1)_q.$ In both cases this entry is $(t_1)_q(t_2)_q$
where $t_1,t_2\in [1,\ldots,l-1]$ and one of $t_1,t_2$ is even
and another is odd. But if $q=\exp (\pi i\frac{l+1}{l})$ and
$t\in [1,\ldots,l-1]$
     then $(t)_q>0 \Leftrightarrow t$ is odd. Finally note that
the entries of $D$ are the entries of $D(j,k)$ divided by some
positive numbers. The Lemma is proved.$\Box$

     \subsubsection{} The above Lemma implies that $\ad_j(k)$
     is a sum of a projective module $P(k)$ and some module
$Y(k,j).$ It follows from ~\ref{resume}
     that $Y(k,j)$ admits a filtration of length 3 with the
following associated graded factors: $L(k); L(-2-k)\oplus
L(2l-2-k); L(k).$

     \subsubsection{} For any $0\le s<l$ an element
$pr_j(K^{-1}E)^s$ is an upper singular vector
     in $\ad_j.$ Let us prove that if $s\le \frac{l-1}{2}$ then
$(X-b_j)pr_j(K^{-1}E)^s\ne 0.$ Consider the action of this
element on $P(J').$ Then $(X-b_j)$ is a surjection onto
$L(J')\subset P(j')$ and the desired result follows from the
fact that $\dim L(J') \ge \frac{l+1}{2}.$ Similarly the vector
$pr_jF^s$ is a lower singular vector
     in $\ad_j,$ and $pr_jF^s\not\in M_j$ if $s\le
\frac{l-1}{2}.$
     For $s=\frac{k}{2}$ we obtain that $\ad_j$ contains a
submodule $L(k)$ which does not lie in $M_j.$ Indeed, the
submodules generated by $pr_j(K^{-1}E)^s$ and $pr_jF^s$
coincide since $\ad_j(k)/M_j(k)\simeq L(k)\oplus L(k)$ and
$\ad_j(k)\supset P(k) \not \subset M_j(k).$

     \subsubsection{} It follows that
     $Y(k,j)=Z(k,j)\oplus L(k)$ where $Z(k,j)$ contains a
submodule $L(k).$ Hence $\ad(k)=\bigoplus \limits _{j\in
H'}\ad_j(k)$ is a direct sum of a few
     copies of $P(k)$ and $Y(k)=\bigoplus \limits _{j\in
H'}Y(k,j)$ where all the subquotients $L(k)$ of $Y(k)$ are the
submodules of $Y(k).$ Now from the autoduality (see
{}~\ref{autoduality}) we see that all the subquotients $L(k)$ are
direct summands of $Y(k).$ Thus $Y(k)$ is a direct sum of its
simple subquotients. Hence in the case ~\ref{resume1} (c) (i.e.
if $k\le 2J$) we have that $$ \ad_j(k)=P(k)\oplus L(k)\oplus
L(k)\oplus L(2l-2-k)\oplus L(-2-k) $$
     This completes the proof of the Main Theorem. $\Box$

     \section{Three matrix lemmas} The results of this section
were used in the previous sections.

     \subsection{} \label{lin1}
     Let $A$ be a $r\times r$ lower-triangular matrix:

     \begin{equation} A=\left( \begin{array}{ccccc} \alpha_1 &
0 & 0 &\ldots & 0 \\ \beta_1 &\alpha_2 & 0 &\ldots &0 \\
\gamma_1 & \beta_2 &\alpha_3 &\ldots & 0 \\ 0 & \gamma_2
&\beta_3 &\ldots & 0 \\ \vdots &\vdots & \vdots &\ddots &\vdots
\\ 0 & 0 &\ldots &\beta_{r-1}&\alpha_r\\ \end{array} \right)
\end{equation}

     {\bf Lemma.} {\em Let $\alpha_i=\alpha_j=\alpha$ for some
$i<j$ and $\alpha_k\ne \alpha$ for $k\ne i,j.$ The matrix $A$
has 2 different eigenvectors with eigenvalue $\alpha$ iff the
determinant of $(j-i)\times (j-i)$ matrix

     \begin{equation} D=\left( \begin{array}{ccccc} \beta_i
&\alpha_{i+1}-\alpha&0 &\ldots & 0 \\ \gamma_i &\beta_{i+1}
&\alpha_{i+2}-\alpha&\ldots &0 \\ 0 & \gamma_{i+1} &\beta_{i+2}
&\ldots & 0 \\ 0 & 0 &\gamma_{i+2} &\ldots & 0 \\ \vdots
&\vdots & \vdots &\ddots &\vdots \\ 0 & 0 &\ldots
&\gamma_{j-2}&\beta_{j-1}\\ \end{array} \right) \end{equation}

     vanishes.}

     {\bf Proof.} Clear.$\Box$

     \subsection{} \label{lin2} Suppose $\alpha_i=\alpha_j$ iff
$i+j=r+1.$ Let $A'$ be the following matrix $$ A'=\left(
\begin{array}{cc}
     A & 0 \\ B' & A \\ \end{array} \right) $$ where $$
B'=\left( \begin{array}{cccc} 0 & 0 &0 & \ldots \\
     1 & 0 &0 & \ldots \\ 0 & 1 & 0 & \ldots \\ \ldots & \ldots
& \ldots & \ldots \\ \end{array} \right) $$

     {\bf Lemma.} {\em Suppose the matrix $A$ above has 2
eigenvectors with eigenvalue $\alpha_i.$ Suppose the matrix $D$
is semisimple of corank 1. Then the matrix $A'-\alpha_i$ has
corank 3.}

     {\bf Proof.} It suffices to consider the case
$i=1,j=r,\alpha_1= \alpha_r=0.$
     Deleting two rows and columns consisting of zeros we
obtain a matrix $$ D'=\left( \begin{array}{cc}
     D & 0 \\ E & D \\ \end{array} \right) $$ where $E$ is a
unit matrix. We have to prove that corank of $D'$ is equal to
1.
     Let $\Image (D)$ (resp. $\Image (D')$) denote the linear
space generated by columns of $D$ (resp. $D'$). Let $pr:\Image
(D')\to {\Bbb C}^{r-1}$ be a map forgetting the last $r-1$
coordinates. Then $pr(\Image (D'))=\Image (D)$ has dimension
$r-2.$ Let us prove that $\Ker (pr)$ has dimension $r-1.$
Indeed, $\Ker (pr)\supset \Image (D)$ and $\Ker (pr)$ contains
the kernel of operator $D.$ Since $D$ is semisimple $\Image
(D)\not \supset \Ker (D).$ The result follows. $\Box$

     \subsection{} \label{lin3}
     {\bf Lemma.} {\em Suppose $D$ is a real matrix, and all
the off-diagonal entries are negative. Then $D$ is semisimple.}

     {\bf Proof.} It is easy to see that conjugating matrix $D$
by some diagonal matrix we can obtain a symmetric matrix.
$\Box$

     \end{document}